\documentclass{amsart}
\usepackage{latexsym,amssymb,amsmath,amsthm,graphics,graphicx,float,psfrag, epsfig, color}
\usepackage{bbm}
\usepackage{amsaddr}

\newcommand{\R}{\mathbb{R}}

\newcommand{\PP}{\mathbb{P}}

\newcommand{\supp}{\mbox{supp}}

\newcommand{\eps}{\varepsilon}
\newcommand{\E}{\mathbb{E}}

\newcommand{\Y}{{Y}}

\newcommand{\HTp}{{H}_{T^\perp}}
\newcommand{\ITp}{{I}_{T^\perp}}
\newcommand{\YTp}{{Y}_{T^\perp}}
\newcommand{\YT}{{Y}_{T}}
\newcommand{\ITop}{{I}_{T_{x_0}^\perp}}
\newcommand{\YTop}{{Y}_{T_{x_0}^\perp}}
\newcommand{\YTo}{{Y}_{T_{x_0}}}

\newcommand{\HT}{{H}_{T}}

\newcommand{\A}{\mathcal{A}}

\newcommand{\Xhat}{\hat{X}}
\newcommand{\xohat}{\hat{x}_0}

\newcommand{\matrixgeq}{\succeq}

\DeclareMathOperator{\tr}{tr}

\newcommand{\aiai}{ a_i a_i^t}
\newcommand{\aipaip}{ a'_i a_i{'^*}}
\newcommand{\fmin}{f_\text{min}}
\newcommand{\ctilde}{\tilde{c}}

\newcommand{\gammatilde}{\tilde{\gamma}}
\newcommand{\xoxot}{x_0 x_0^t}
\newcommand{\xxt}{x x^t}
\newcommand{\ASc}{\A_{S^c}}
\newcommand{\AS}{\A_S}
\newcommand{\Tp}{T^\perp}
\DeclareMathOperator{\sgn}{sgn}
\DeclareMathOperator{\polylog}{polylog}
\newcommand{\Sc}{S^c}
\newcommand{\etaSc}{\eta_{\Sc}}
\newcommand{\etaS}{\eta_S}
\newcommand{\Splus}{S^+}
\newcommand{\Sminus}{S^-}
\newcommand{\indic}{\mathbbm{1}}
\newcommand{\YplusTp}{\Y_{+, \Tp}}
\newcommand{\YplusminusTp}{\Y_{\pm, \Tp}}
\newcommand{\YminusTp}{\Y_{-, \Tp}}
\newcommand{\YoTp}{\Y_{0, \Tp}}
\newcommand{\YplusT}{\Y_{+, T}}
\newcommand{\YplusminusT}{\Y_{\pm, T}}

\newcommand{\YoT}{\Y_{0, T}}
\newcommand{\Yplus}{\Y_+}
\newcommand{\Yminus}{\Y_-}
\newcommand{\Yo}{\Y_0}
\newcommand{\Etilde}{\tilde{E}}
\newcommand{\Neps}{\mathcal{N}_\eps}

\newtheorem{theorem}{Theorem}
\newtheorem{lemma}{Lemma}

\title{PhaseLift is robust to a constant fraction of arbitrary errors}
\author{Paul Hand }
\address{Rice University\\ Department of Computational and Applied Mathematics, \\ MS-134, 6100 Main Street, Houston, TX 77005}
\email{hand@rice.edu}
\date{February 14, 2015}		

\begin{document}
\maketitle

\begin{abstract}
Consider the task of recovering an unknown $n$-vector from phaseless linear measurements.  This task is the phase retrieval problem.  Through the technique of lifting, this nonconvex problem may be convexified into a semidefinite rank-one matrix recovery problem, known as PhaseLift.  Under a linear number of exact Gaussian measurements, PhaseLift recovers the unknown vector exactly with high probability.  Under noisy measurements, the solution to a variant of PhaseLift has error proportional to the $\ell_1$ norm of the noise.  In the present paper, we study the robustness of this variant of PhaseLift to a case with noise and gross, arbitrary corruptions.  We prove that PhaseLift can tolerate a  small, fixed fraction of gross errors, even in the highly underdetermined regime where there are only $O(n)$ measurements.  The lifted phase retrieval problem can be viewed as a rank-one robust Principal Component Analysis (PCA) problem under generic rank-one measurements.  From this perspective, the proposed convex program is simpler that the semidefinite version of the sparse-plus-low-rank formulation standard in the robust PCA literature.  Specifically, the rank penalization through a trace term is unnecessary, and the resulting optimization program has no parameters that need to be chosen.  The present work also achieves the information theoretically optimal scaling of $O(n)$ measurements without the additional logarithmic factors that appear in existing general robust PCA results.
\end{abstract}

\smallskip
\noindent \textbf{Acknowledgements.} 
The author is partially supported by generous funding from the National Science Foundation.


\section{Introduction}
This paper establishes robustness of an algorithm for recovering a vector $x_0 \in \R^n$ from phaseless linear measurements that contain noise and a constant fraction of gross, arbitrary errors.  That is, for fixed measurement vectors $a_i \in \R^n$ for $i = 1 \ldots m$, our task is to find $x_0$ satisfying 
\begin{align}
b_i = | \langle x_0, a_i \rangle|^2 + \eps_i  + \eta_i \label{measurements-components}
\end{align}
for known $b_i \in \R$, known $a_i$, and unknown $\eta_i$ and $\eps_i$.  Here $\eta_i$ will represent the noise in the measurements, and $\eps_i$ will represent gross, arbitrary errors.  This recovery problem is known as  phase retrieval.  Measurements of form \eqref{measurements-components} arise in several applications,  such as X-ray crystallography, optics, and microscopy \cite{harrison1993phase, millane1990phase,CESV2011}.  In such applications, extremely large errors in some measurements may be due to sensor failure, occlusions, or other effects.  Ideally, recovery algorithms could provably tolerate a small number of such errors.

Recently, researchers have introduced algorithms for the phase retrieval problem that have provable recovery guarantees \cite{CESV2011,CSV2013}.  The insight of these methods is that the phase retrieval problem can be convexified by lifting it to the space of matrices.  That is, instead of searching for the vector $x_0$, one can search for the lifted matrix $\xoxot$.  The quadratic measurements \eqref{measurements-components} then become linear measurements on this lifted matrix.  As the desired matrix is semidefinite and rank-one, one can write a rank minimization problem under the semidefinite and data constraints, which has a convex relaxation known as PhaseLift.  In this noiseless case, PhaseLift is the program
\begin{align}
\min_X \text{tr}(X) \text{ subject to } X \matrixgeq 0, \{a_i^t X a_i = b_i\}_{i = 1 \ldots m} \label{phaselift-trace}
\end{align}
Here, the trace of $X$ is a convex proxy for the rank of a positive semidefinite $X$.  An estimate for the underlying signal $x_0$ can be computed by the leading eigenvector of the optimizer of \eqref{phaselift-trace}.

As in \cite{CSV2013, DH2014, CL2012}, we will seek recovery guarantees for independent identically distributed Gaussians
$$
a_i \sim \mathcal{N}(0, I_n).
$$
Under this data model,  \cite{DH2014} and \cite{CL2012} have shown that \eqref{phaselift-trace}  can be simplified to the semidefinite feasibility problem
\begin{align}
\text{find } X \matrixgeq 0 \text{ such that } \{a_i^t X a_i = b_i\}_{i = 1 \ldots m} \notag 
\end{align}
This feasibility problem succeeds at finding $\xoxot$ exactly  with high probability when $m \geq c n$ for a sufficiently large $c$  \cite{CL2012}.  This scaling is quite surprising because there are only $O(n)$ measurements for an $O(n^2)$ dimensional object.  As discussed in \cite{DH2014}, the semidefinite cone is sufficiently `pointy' that the high-dimensional affine space of data-consistent matrices intersects the semidefinite cone only at exactly one point.

In the noisy case  without gross erros, that is for $\eps=0$, \cite{CL2012} showed that the PhaseLift variant
\begin{align}
\min \sum_{i} | a_i^t X a_i - b_i | \text{ subject to } X \matrixgeq 0 \label{phaselift-l1} 
\end{align}
successfully recovers a matrix near $\xoxot$ with high probability.  Specifically, they prove that the solution $\Xhat$ to \eqref{phaselift-l1} satisfies $\|\Xhat - \xoxot\|_\text{F} \leq C_0 \|\eta \|_1/m$ with high probability.  From $\Xhat$, an estimate of $x_0$ can be obtained by $$\xohat = \sqrt{\hat{\lambda}_1} \hat{u}_1$$
where $(\hat{\lambda}_1,  \hat{u}_1)$ is the leading eigenvector and eigenvalue pair for $\Xhat$.   In \cite{CL2012}, the authors prove that $|\xohat - \pm x_0 | \leq C_0 \min(\|x_0\|, \|\eta\|_1 / m \|x_0\| )$ for some $C_0$.  

The contribution of the present paper is to show that the program \eqref{phaselift-l1} is additionally robust against a constant fraction of arbitrary errors.    For a fixed set of coefficients that contain gross errors, we show that approximate recovery succeeds with high probability for arbitrary signals and arbitrary values of the gross errors.   
\begin{theorem} \label{theorem-exact-recovery}
There exist positive numbers $\fmin, \gamma, c, C, C'$ such that the following holds.  Let $m \geq cn$.  Fix a set $S \subset \{1 \ldots m\}$ such that $|S|/m \leq \fmin$.  On an event of probability at least $1 - e^{-\gamma m}$, for any $x_0 \in \R^n$ and for any $\eps$ with $\supp(\eps) \subseteq S$, the minimizer $\Xhat$ to \eqref{phaselift-l1} satisfies $$\|\Xhat - \xoxot\|_\text{F} \leq C \frac{\|\eta\|_1}{m}.$$
The resulting estimate for $x_0$ satisfies
$$\| \xohat - \pm x_0 \| \leq C' \min \left( \|x_0\|, \frac{\| \eta\|_1}{m \|x_0\|}  \right)$$

\end{theorem}
Note that this high-probability result is universal over $x_0$ and $\eps$ and does not only apply for merely for a fixed signal or for a fixed error vector $\eps$.

In the case of no gross errors, in which $\eps = 0$, this theorem reduces to the result in \cite{CL2012} mentioned above.  In the noiseless case, in which $\eta = 0$, the theorem guarantees exact recovery of $x_0$ with high probability under a a linear number of measurements, of which a constant fraction are corrupted.  

We now explore the optimality of this theorem.  The scaling of $m$ versus $n$ is information theoretically optimal and has no unnecessary logarithmic factors.  The noise scaling is the same as in \cite{CL2012}, and its optimality was established there.  For arbitrary errors, the fixed fraction of gross errors can not be extended to a case where $\fmin \geq 1/2$ because one could build a problem where half of the measurements are due to an $x_0$ and the other half are due to some $x_1$.  In such a case, recovery would be impossible.

\subsection{Relation to Robust PCA}

Much recent work in matrix completion has studied the recovery of low-rank matrices from arbitrary corruptions to its entries, known as robust Principal Component Analysis (PCA).  Results in this framework typically involve measuring some of the entries of a low rank $n \times n$ matrix $X$ and assuming that some fraction of those measurements are arbitrarily corrupted, giving the data matrix $A$.  The matrix $X$ can then be recovered under certain conditions by a sparse-plus-low-rank convex program:
\begin{align}
\min \lambda \|X\|_* + \|E\|_1 \text{ such that } \mathcal{P}(X + E) = \mathcal{P}(A) \label{sparse-plus-low-rank}
\end{align}
where $\|X\|_*$ is the nuclear norm of $X$, $\lambda$ is a constant, $\|E\|_1$ is the $\ell_1$ norm of the vectorization of $E$, and $\mathcal{P}$ is the projection of a matrix onto the observed entries. \cite{CLMW2011, CSPW2011, HKZ2011, WGRPM2009, GWLCM2010, CJSC2013, Li2013}.  Results from this formulation have been quite surprising.  For example, under an appropriate choice of $\lambda$ and under an incoherence assumption, the sparse-plus-low-rank decomposition succeeds for sufficiently low rank $X$ when $O(n^2)$ entries are measured and a small fraction of them have arbitrary errors \cite{CLMW2011}.  Subsequent results have been proved that only require $m \gtrsim rn \polylog(n)$ measurements, where $r$ is the rank of $X$ \cite{CJSC2013, Li2013}.  This result is information theoretically optimal except for the polylogarithmic factor.  

The present paper can be viewed as a rank-one semidefinite robust PCA problem under generic rank-one measurements.  From this perspective, we would naturally formulate the PhaseLift problem under gross errors by
\begin{align}
\min \lambda \tr(X) +  \sum_{i} | a_i^t X a_i - b_i | \text{ subject to } X\matrixgeq 0. \label{phaselift-l1-trace}
\end{align}
The present paper shows that explicit rank penalization by the trace term is not fundamental for exact recovery in the presence of arbitrary errors.  That is, \eqref{phaselift-l1-trace} can be simplified by taking $\lambda = 0$.  The resulting program has no free parameters that need to be explicitly tuned.  As in \cite{DH2014, CL2012}, the positive semidefinite cone provides enough of a constraint to enforce low-rankness.  The present paper also shows that rank-one matrix completion can succeed under an information theoretically optimal data scaling.  Specifically, the extra logarithmic factors from low-rank matrix completion and robust PCA do not appear in Theorem \ref{theorem-exact-recovery}.  The present work also differs from the standard robust PCA literature in that the measurements are generic and are not direct samples of the entries of the unknown matrix.  

\subsection{Numerical simulation}

\begin{figure}[h]
	{\large
%
%
\begin{psfrags}%
\psfragscanon%
%
\psfrag{s03}[t][t]{\color[rgb]{0,0,0}\setlength{\tabcolsep}{0pt}\begin{tabular}{c}Number of measurements ($m$)\end{tabular}}%
\psfrag{s04}[b][b]{\color[rgb]{0,0,0}\setlength{\tabcolsep}{0pt}\begin{tabular}{c}Signal length ($n$)\end{tabular}}%
\psfrag{s05}[l][l]{\color[rgb]{0,0,0}\setlength{\tabcolsep}{0pt}\begin{tabular}{l}{\Large $m = \frac{(n+1)n}{2}$}\end{tabular}}%
%
\psfrag{x01}[t][t]{50}%
\psfrag{x02}[t][t]{100}%
\psfrag{x03}[t][t]{150}%
\psfrag{x04}[t][t]{200}%
\psfrag{x05}[t][t]{250}%
%
\psfrag{v01}[r][r]{50}%
\psfrag{v02}[r][r]{25}%
%
\resizebox{8cm}{!}{\includegraphics{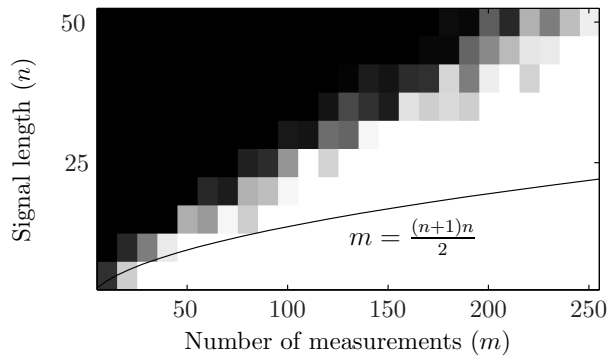}}%
\end{psfrags}%
%

  	}
  \caption{Recovery error for the PhaseLift matrix recovery problem \eqref{phaselift-l1} as a function of $n$ and $m$, when 5\% of measurements contain large errors.  Black represents an average recovery error of 100\%.  White represents zero average recovery error.  Each block corresponds to the average from 10 independent trials.  The solid curve depicts when the number of measurements equals the number of degrees of freedom in a symmetric $n \times n$ matrix.    The number of measurements required for successful recovery appears to be linear in $n$, even with a small fraction of large errors. }
    \label{figure-phase-transition}
\end{figure}

We now explore the empirical performance of \eqref{phaselift-l1} by numerical simulation.  Let the signal length $n$ vary from 5 to 50, and let the number of measurements $m$ vary from 10 to 250.  Let $x_0 = e_1$.  For each $(n,m)$, we consider measurements such that
\begin{align}
\begin{cases} b_i \sim  \text{Uniform}([0, 10^4]) & \text{ if } 1 \leq i \leq \lceil 0.05 m \rceil, \\ b_i =  |\langle x_0, a_i \rangle |^2  & \text{ otherwise.} \end{cases}
\end{align}
We attempt to recover $\xoxot$ by solving \eqref{phaselift-l1} using the SDPT3  solver \cite{TTT1999, TTT2003} and YALMIP \cite{L2004}.  For a given optimizer $\Xhat$, define the capped relative error as $$\min(\|\Xhat - \xoxot\|_\text{F} / \|\xoxot\|_\text{F}, 1).$$  Figure \ref{figure-phase-transition} plots the average capped relative error over 10 independent trials.  It provides empirical evidence that the matrix recovery problem \eqref{phaselift-l1} succeeds under a linear number of measurements, even when a constant fraction of them contain very large errors.  

\section{Proofs}

Let $\A: \mathcal{S}^n \to \R^m$ be defined by the mapping $X \mapsto (a_i^t X a_i)_{i = 1 \ldots m}$, where $\mathcal{S}^n$ is the space of symmetric real-valued $n \times n$ matrices.  Note that $\A^* \lambda = \sum_i \lambda_i \aiai$.  Let $\AS$ be the restriction of $\A$ onto the coefficients given by the set $S$. Let $e_i$ be the $i$th standard basis vector.  Let $X_0 = x_0 x_0^t$.  We can write the measurements \eqref{measurements-components} as $$b = \A X_0 + \eps + \eta.$$  Similarly, the optimization program \eqref{phaselift-l1} can be written as $$ \min \|\A X  - b \|_1 \text{ such that } X \matrixgeq 0.$$

We introduce the following notation.  Let $\|X\|_1$ be the nuclear norm of the matrix $X$.  When $X \matrixgeq 0$, $\|X\|_1 = \tr(X)$.  Denote the Frobenius and spectral norms of $X$ as $\|X\|_\text{F}$ and $\|X\|$, respectively.  Given $x_0$, let $T_{x_0} = \{ y x_0^t + x_0 y^t \mid y \in \R^n \}$.  Note that $T_{e_1}$ is the space of symmetric matrices supported on their first row and column.  The orthogonal complement $T_{e_1}^\perp$ is then the space of matrices supported in the lower-right $n-1 \times n-1 $ block.  When $x_0$ is clear, we will simply write $T$ instead of $T_{x_0}$.  Let $I$ be the identity matrix, and let $\indic(E)$ be the indicator function of the event $E$. 

\subsection{Recovery by dual certificates}

The proof of Theorem \ref{theorem-exact-recovery} will be based on dual certificates, as in \cite{CSV2013, DH2014, CL2012}.  A dual certificate is an optimal variable for the dual problem to \eqref{phaselift-l1}.  Its existence certifies the correctness of a minimizer to \eqref{phaselift-l1}.

The first order optimality conditions at $X_0$ for \eqref{phaselift-l1} are given by
\begin{align}
\tilde{Y} &= \A^* \tilde{\lambda} \\
\tilde{\lambda} &\in \partial \|\cdot \|_1(-\eps)\\
\tilde{Y} &\matrixgeq 0 \label{psd-constraint}\\
\langle \tilde{Y}, X_0 \rangle &= 0 \label{slackness-constraint}
\end{align}
where $\partial \|\cdot \|_1(-\eps)$ is the subgradient of the $\ell_1$ norm evaluated at $-\eps$.  Note that \eqref{psd-constraint} and \eqref{slackness-constraint} imply $\tilde{Y}_T = 0$.   Such a $\tilde{Y}$ would be dual certificate for \eqref{phaselift-l1}.  Unfortunately, constructing such a $\tilde{Y}$ that exactly satisfies these conditions is difficult.  As in \cite{CSV2013, DH2014, CL2012}, we seek an inexact dual certificate, which approximately satisfies these conditions.  Specifically, we will build a dual certificate $Y = \A^* \lambda$ that satisfies
\begin{align}
&\YTp \matrixgeq \ITp \label{ytp-condition}\\
&\| \YT\|_\text{F} \leq 1/2 \label{yt-condition} \\
&\begin{cases}\lambda_i  = - \frac{7}{m} \sgn(\eps_i)  & \text{ if } \eps_i \neq 0 \\ 
\lambda_i \leq \frac{7}{m} & \text{ if } \eps_i = 0.
\end{cases} \label{lambda-condition}
\end{align}

To prove that existence of such a $Y$ will guarantee successful recovery of $\xoxot$ with high probability, we will rely on two technical lemmas.  The first technical lemma provides $\ell_1$-isometry bounds on $\A$ and was proven in \cite{CSV2013}.
\begin{lemma}[\cite{CSV2013}]\label{lemma-l1-isometry}
There exist constants $c_0, \gamma_0$ such that if $m \geq c_0 n$, then with probability at least $1 - e^{-\gamma_0 m}$,
\begin{align}
\frac{1}{m} \|\A(X) \|_1 &\leq \left(1 + \frac{1}{16} \right) \| X \|_1 \text{ for all } X, \\
\frac{1}{m} \|\A(X) \|_1 &\geq 0.94 \left(1 - \frac{1}{16}\right) \|X\| \text{ for all } \text{symmetric, rank-2 } X
\end{align}
\end{lemma}

We will need simultaneous control of the $\ell_1$-isometry properties over several subsets of measurements.

\begin{lemma}\label{lemma-l1-isometry-subsets}
There exists a constant $\gammatilde_0$ such that the following holds.  Let $m \geq 100 c_0 n$, and fix a support set $S$ with $|S| = \lceil 0.01 m \rceil$.  There is an event $E_S$ with probability at least $1 - e^{-\gammatilde_0 m}$ on which
\begin{align}
\frac{1}{|S^c|} \| \ASc X \|_1 &\geq 0.94 \left( 1 - \frac{1}{16} \right) \|X\| \text{ for all symmetric rank-2 } X,\\
\frac{1}{|S|} \| \AS X\|_1 &\leq \left( 1 + \frac{1}{16} \right) \| X\|_1 \text{ for all } X,\\
\frac{1}{|m|} \| \A X\|_1 &\leq \left( 1 + \frac{1}{16} \right) \| X\|_1 \text{ for all } X.
\end{align}

\end{lemma}
The proof of Lemma \ref{lemma-l1-isometry-subsets} is immediate from Lemma \ref{lemma-l1-isometry}.

In order to prove that a dual certificate guarantees recovery, we establish a technical result that an optimal solution $X_0+H$ to \eqref{phaselift-l1} lies near the cone $\|\HTp\|_1 \leq \frac{1}{2} \| \HT\|_\text{F}$ with high probability.  This property is a strong version of injectivity on $T$.

\begin{lemma} \label{lemma-feasible-cone}
Fix a support set $S$ with $|S| = \lceil 0.01 m \rceil$  and let $m \geq 100 c_0 n$.  On the event $E_S$ from Lemma \ref{lemma-l1-isometry-subsets}, for all $x$ and for all $\eps$ with supp$(\eps) \subseteq S$, any optimal $X_0 + H$ satisfies $\|\HTp \|_1 \geq 0.56 \| \HT \|_F - \frac{2}{m}\|\eta\|_1$.
\end{lemma}

\begin{proof}
By assumption, $\| \A H - \eps - \eta \|_1 \leq \|\eps + \eta\|_1$.    By the additivity of the $\ell_1$ norm over vectors with disjoint supports,
\begin{align}
\|\ASc H - \etaSc \|_1 + \|\AS H - \eps - \etaS\|_1 \leq \|\eps + \etaS\|_1 + \|\etaSc\|_1.
\end{align}
By the triangle inequality, we have
\begin{align}
&\| \ASc H  - \etaSc \|_1 + \|\eps + \etaS\|_1 - \| \AS H \|_1 \leq \| \eps + \etaS \|_1 + \| \etaSc\|_1,
\end{align}
which implies
\begin{align}
&\| \ASc H\|_1 \leq 2 \| \etaSc\|_1 + \| \AS H \|_1 \label{asch-vs-ash}
\end{align}

Breaking \eqref{asch-vs-ash} into its components on $T$ and $\Tp$ and applying the triangle inequality, we have
\begin{align}
&\|\ASc \HT \|_1 - \| \ASc \HTp \|_1  \leq 2 \| \etaSc\|_1 +  \| \AS \HT \|_1 + \| \AS \HTp\|_1\\
\Rightarrow &\|\ASc \HT\|_1 \leq 2 \| \etaSc\|_1 +  \| \AS \HT\|_1 + \| \A \HTp \|_1 \label{ascht-vs-asht-ahtp}
\end{align}
We now apply the $\ell_1$ isometry bounds from Lemma \ref{lemma-l1-isometry} on each term of \eqref{ascht-vs-asht-ahtp}.  On the event $E_S$, 
\begin{align}
\| \ASc \HT \|_1 &\geq 0.94 \left(1 - \frac{1}{16} \right) |S^c| \|\HT\| \\
&\geq 0.94 \left(1 - \frac{1}{16}\right) \frac{|S^c|}{\sqrt{2}} \| \HT\|_\text{F}, \label{ascht-bound}
\end{align}
where the second inequality follows because $\HT$ has rank at most 2.
On the event $E_S$, 
\begin{align}
\|\AS \HT \|_1 &\leq |S| \left(1 + \frac{1}{16}\right) \|\HT\|_1 \\
&\leq |S| \left(1 + \frac{1}{16} \right) \sqrt{2} \|\HT\|_\text{F}
\end{align}
On the event $E_S$,
\begin{align}
\|\A \HTp\|_1 \leq m \left(1 + \frac{1}{16}\right) \| \HTp \|_1 \label{ahtp-bound}
\end{align}
Combining \eqref{ascht-vs-asht-ahtp}--\eqref{ahtp-bound}, we have
\begin{align}
\left(\frac{0.94}{\sqrt{2}} \frac{\left(1 - \frac{1}{16} \right)}{\left( 1 + \frac{1}{16} \right)} \frac{|S^c|}{m}  - \sqrt{2} \frac{|S|}{m} \right) \|\HT\|_\text{F} \leq \frac{2}{m} \|\etaSc\|_1 + \| \HTp \|_1
\end{align}
Thus, $0.56 \| \HT \|_F \leq \frac{2}{m} \|\eta \|_1 + \| \HTp\|_1$ on the event $E_S$. 
\end{proof}

We may now prove that existence of an inexact dual certificate \eqref{ytp-condition}--\eqref{lambda-condition} will guarantee successful recovery of a matrix near $X_0 = \xoxot$ with high probability, provided that there are few enough arbitrary errors.
\begin{lemma} \label{lemma-recovery-with-dual-certificate}
There exists a $C$ such that the following holds.  Fix $S$ such that $|S| = \lceil 0.01 m\rceil$.  Let $m \geq 100 c_0 n$. Fix $x_0 \in \R^n$.   Fix $\eps \in \R^m$ such that $\supp(\eps) \subseteq S$.  Suppose that there exists $Y = \A^* \lambda$ satisfying \eqref{ytp-condition}--\eqref{lambda-condition}.  Then, on the event $E_S$ from Lemma \ref{lemma-l1-isometry-subsets},  a minimizer $\Xhat$ of  \eqref{phaselift-l1} satisfies $\|\Xhat - X_0\|_\text{F} \leq C\frac{\|\eta\|_1}{m}$.
\end{lemma}
\begin{proof}
Let $\Xhat = X_0 + H$ be a minimizer for \eqref{phaselift-l1}, which implies  that $\|\A(X_0 + H) - b \|_1 \leq \| \A(X_0) - b\|_1$.  That is $$ \|\A H -  \eps - \eta\|_1 \leq \|\eps + \eta\|_1.$$
Letting $\alpha = 7/m$, condition \eqref{lambda-condition} gives
\begin{align}
\lambda/\alpha \in \partial \|\cdot\|_1(-\eps)
\end{align}
Hence,
\begin{alignat}{2}
&& \|-\eps\|_1 + \langle \lambda/\alpha, \A H - \eta \rangle &\leq \|\eps\|_1 + \|\eta\|_1 \\
&\Rightarrow  &\langle \lambda, \A H - \eta \rangle &\leq \alpha \|\eta\|_1 \\
&\Rightarrow  &\langle Y, H \rangle &\leq \langle \lambda, \eta \rangle + \alpha \|\eta\|_1 \\
&\Rightarrow  &\langle Y, H \rangle &\leq 2 \alpha \|\eta\|_1 \label{yh-negative}
\end{alignat}
Decomposing \eqref{yh-negative} into $T$ and $\Tp$, we have
\begin{align}
\langle \YT, \HT \rangle \leq - \langle \YTp, \HTp \rangle + 2 \alpha \|\eta\|_1
\end{align}
As $\YTp \matrixgeq 0$ and $\HTp \matrixgeq 0$, we have
\begin{align}
\langle \YTp, \HTp \rangle \leq  | \langle \YT, \HT \rangle | + 2 \alpha \|\eta\|_1
\end{align}
By conditions \eqref{ytp-condition}--\eqref{yt-condition}
\begin{align}
\|\HTp \|_1 \leq | \langle \YTp, \HTp\rangle | \leq | \langle \YT, \HT \rangle | + 2 \alpha \|\eta\|_1 \leq \frac{1}{2} \|\HT\|_F + 2 \alpha \|\eta\|_1 \label{cone-condition-first}
\end{align}
By Lemma \ref{lemma-feasible-cone}, on the event $E_S$, \begin{align}
\|\HTp \|_1 \geq 0.56 \| \HT \|_\text{F} - \frac{2}{m} \|\eta\|_1. \label{cone-condition}
\end{align}
Combining \eqref{cone-condition-first} and \eqref{cone-condition}, and using $\alpha = 7/m$, we get
\begin{align}
0.56 \|\HT\|_\text{F} \leq \frac{2}{m} \| \eta\|_1 + 0.5 \| \HT\|_\text{F} + \frac{14}{m} \| \eta\|_1
\end{align}
So,
\begin{align}
\|\HT\|_\text{F} \leq \frac{C'}{m} \| \eta\|_1 \text{ and thus } \|\HTp\|_\text{F} \leq \frac{C''}{m} \| \eta\|_1
\end{align}
We conclude $\|\Xhat - X_0\|_\text{F} = \|H\|_\text{F} \leq \frac{C}{m} \|\eta\|_1$ for some $C$.
\end{proof}

\subsection{Construction of the dual certificate}

We now construct the dual certificate for arbitrary $x_0$.  Our construction will be a modification to the dual certificate in \cite{CL2012}.  Also similar to \cite{CL2012}, we will build dual certificates with high probability on a net of $x_0$.  We will then use a continuity argument to get a dual certificate for a arbitrary $x_0$.  

Let $\Splus$ and $\Sminus$ be disjoint supersets of the indices over which $\eps$ is positive or negative, respectively.  Let $S = \Splus \cup \Sminus$.  For pedagogical purposes, $\Splus$ and $\Sminus$ should be thought of as exactly the indices over which $\eps$ is positive or negative.  For technical reasons, we let them be supersets of cardinality linear in $n$, in order to use standard probability bounds. 
For a fixed choice of $\Splus$ and $\Sminus$, let the inexact dual certificate $Y$ be defined by
\begin{align}
Y = \frac{1}{m} \Bigl[ \sum_{i \in \Splus} -7 \aiai  + \sum_{i \in \Sminus} 7 \aiai  + \sum_{i \in \Sc} [\beta_0 - |\langle a_i, \frac{x_0}{\|x_0\|}\rangle|^2 \indic(|\langle a_i, \frac{x_0}{\|x_0\|} \rangle | \leq 3 )] \aiai    \Bigr] \label{dual-certificate}
\end{align}
where $\beta_0 = \mathbb{E} z^4 \mathbbm{1}(|z| \leq 3) \approx 2.6728$ and $z$ is a standard normal random variable.  We will refer to each of the terms of in the right hand side of \eqref{dual-certificate} as $\Yplus, \Yminus$, and $\Yo$, respectively.  

The form of $\Yo$ is due to \cite{CL2012}, and the intuition behind it is as follows.  Note that $\E (\aiai) = I_n$ and $\E ( |\langle a_i, e_1 \rangle|^2 \aiai) = \begin{pmatrix}\tilde{\beta_0} & 0 \\ 0 & I_{n-1} \end{pmatrix}$, where $\tilde{\beta_0} = \E z^4$ for a standard normal $z$.  The construction $\frac{1}{m} \sum_{i=1}^m [\tilde{\beta}_0 - |\langle a_i, e_1\rangle|^2 ] \aiai $ thus has expected value $(\tilde{\beta_0} - 1) \ITp$, which would provide the exact dual certificate conditions \eqref{psd-constraint}--\eqref{slackness-constraint}.  As shown in \cite{CL2012}, a satisfactory inexact dual certificate can be built with $m = O(n)$ and coefficients that are truncated to be no larger than $7/m$.  In the present formulation, the terms $\Yplus$ and $\Yminus$ are then set to have coefficients $\mp 7/m$ in order to satisfy  \eqref{lambda-condition}.

We now show that for a fixed signal and a fixed pattern of signs of $\eps$, that a dual certificate exists with high probability.  

\begin{lemma} \label{lemma-dual-certificate-construction}
There exists constants $c, \gamma^*$ such that the following holds.  Let $m \geq c n$.  Fix $x_0 \in \R^n$ and $\eps \in \R^m$.   Let $\Splus$ and $\Sminus$ be fixed disjoint sets of cardinality $\lceil 0.001 m \rceil$.  Then the dual certificate $Y$ from \eqref{dual-certificate} satisfies \eqref{lambda-condition}, $\|\YT\|_\text{F} \leq 1/4$, and $\|\YTp - \frac{17}{10} \ITp \| \leq \frac{3}{10}$ with probability at least $1 - e^{-\gamma^* m}$. 
\end{lemma}

\begin{proof}
Without loss of generality, it suffices to assume $x_0 = e_1$.
It suffices to show that with high probability
\begin{align}
\left \| \YoTp - \frac{17}{10} \ITp \right \| &\leq 0.15, \label{YoTp-bound} \\
\| \YoT\|_\text{F} &\leq \frac{3}{20}, \label{YoT-bound}\\
\| \YplusminusTp\| &\leq 0.015, \label{YplusTp-bound}\\
\| \YplusminusT\|_\text{F} &\leq 0.035. \label{YplusT-bound}
\end{align}

First, we establish \eqref{YoTp-bound}--\eqref{YoT-bound}.  By Lemma 2.3 in \cite{CL2012}, there exist $\ctilde, \gammatilde$ such that if $|S^c| \geq \ctilde n$, then with probability at least $1-e^{-\gammatilde |S^c|},$
\begin{align}
\left \| \frac{m}{|S^c|} \YoTp - \frac{17}{10} \ITp \right \| &\leq 1/10, \text{ and }\\
\left \| \frac{m}{|S^c|} \YoT \right \| &\leq 3/20. \label{yot-bound-two}
\end{align}
Thus,
\begin{align}
\left \|\YoTp - \frac{17}{10} \ITp \right \| &\leq \left \| \YoTp - \frac{|S^c|}{m} \frac{17}{10} \ITp \right \|  \notag  + \left( 1 - \frac{|S^c|}{m} \right) \frac{17}{10} \leq  0.15
\end{align}
which establishes \eqref{YoTp-bound}.  By \eqref{yot-bound-two}, we get \eqref{YoT-bound} immediately.

Next, we establish \eqref{YplusTp-bound}.  Let $a'$ be the vector formed by the last $n-1$ components of $a$.  Observe that $\sum_{i \in \Splus} \aipaip$ is a Wishart random matrix.  Standard estimates for singular values of random matrices with Gaussian i.i.d. entries, such as Corollary 5.35 in \cite{V2012}, apply.  If $|\Splus| = \lceil 0.001 m \rceil \geq \ctilde_0 n$, with probability at least $1 - e^{-\gammatilde_1 m}$ for some $\gammatilde_1$,
$$
\left \| \frac{1}{|\Splus|} \sum_{i \in \Splus} \aipaip - \ITp \right \| \leq 1/2 
$$
Hence,
$$
\left \| \frac{1}{|\Splus|} \sum_{i \in \Splus} \aipaip  \right \| \leq 3/2
$$
Thus,
$
 \| \YplusTp \| \leq \frac{3}{2} \frac{|\Splus|}{m},
$
and we arrive at \eqref{YplusTp-bound}.  The bound for the $\YminusTp$ term is identical.  

Next, we establish \eqref{YplusT-bound}. Note that $\Yplus = -7 \frac{|\Splus|}{m} \cdot \frac{1}{|\Splus|} \sum_{i \in \Splus} \aiai$.  As per Lemma \ref{lemma-aiai-bounds}, if $|\Splus| = \lceil 0.001 m \rceil \geq \ctilde_1 n$,  then $\| \YplusT\|_\text{F} \leq 7 \frac{|\Splus|}{m} \cdot 5 \leq 7 \cdot 0.001 \cdot 5$, with probability at least $1 - e^{-\gammatilde_1 m}$.  

Thus \eqref{YoTp-bound}--\eqref{YplusT-bound} hold simultaneously with probability at least $1 - e^{-\gamma^* m}$ for some $\gamma^*$ provided that $c \geq {\max(2 \ctilde, 1000 \ctilde_0, 1000 \ctilde_1)}.$

\end{proof}

The behavior of $\YplusminusT$ relies on the following probability estimate for the behavior of a Gaussian Wishart matrix on $T$.

\begin{lemma} \label{lemma-aiai-bounds}
Let $x_0 = e_1$.    Let $A = \frac{1}{m} \sum_{i=1}^m \aiai$. There exists $\ctilde_1, \gammatilde_1$ such that if $m \geq \ctilde_1 n$ then $\| A_T \|_\text{F} \leq 5$ with probability at least $1 - e^{-\gammatilde_1 m}$.  

\end{lemma}

\begin{proof}
Let $y$ be the $(1,1)$ entry of $A$, and let $y'$ be the rest of the first column of $A$.  Then $\|A_T\|_\text{F}^2 = y^2 + \|y'\|_2^2.$  So, $y = \frac{1}{m} \sum_{i=1}^m a_i^2(1)$.  Hence, $my \sim \chi^2_m$.  Standard results on the concentration of chi-squared variables, such as Lemma 1 in  \cite{LM2000}, give 
$$
\PP(my \geq 4 m) \leq e^{-\gammatilde_{1,2} m}.
$$
for some $\gammatilde_{1,2}$.  Hence $\PP(y^2 \geq 16) \leq e^{-\gammatilde_{1,2} m}$.

Now, it remains to bound $\|y'\|_2^2$. We can write $y' = \frac{1}{m} Z' c$, where $Z' = [a'_1, \ldots, a'_m]$ and $c_i = a_i(i)$, where $Z'$ and $c$ are independent.  Note that $\|c\|_2^2$ is a Chi-squared random variable with $m$ degrees of freedom.  Hence, with probability at least $1-e^{- \gammatilde_{1, 2} m}$,
$$
\|c\|_2^2 \leq 4m
$$
For fixed $\|x\|_2=1$, $\|Z'x\|_2^2 \sim \chi^2_{n-1}$, and hence with probability at least $1 - e^{-\gamma_{1,3} m}$
$$
\|Z' x \|_2^2 \leq m
$$
when $m \geq \ctilde_1 n$. Hence, $ m^2 \| y' \|_2^2 \leq 4 m \cdot m$
with probability at least $1 -2 e^{-\gammatilde_{1,4} m}$. So, $\|y'\|_2 < 2$ with probability at least $1 -2 e^{-\gammatilde_{1,4} m}$.

So $\|A_T\|_\text{F}^2 \leq 25$, and hence $\|A_T\|_\text{F} \leq 5$ with probability at least $1 - e^{-\gammatilde_1 m}$ for some $\gammatilde_1$. 

\end{proof}

We now show that for a fixed signal and support set of gross errors, there is a high probability that a dual certificate exists simultaneously for all gross errors.

\begin{lemma} \label{dual-certificate-over-all-errors}
Fix $x_0$ and a support set $S$.  If $m \geq cn$ and $|S|/m \leq \min(0.001, \gamma^*/ 2 \log 2)$, then there is an event $\Etilde_{S,x_0}$ on which for all $\eps$ with $\supp(\eps) \subseteq S$, there exists a $Y$ satisfying \eqref{lambda-condition}, $\|\YT\|_\text{F} \leq 1/4$, and $\| \YTp - \frac{17}{10} \ITp \| \leq 3/10$.  The probability of $\Etilde_{S,x_0}$ is at least $1 - e^{-\gamma^* m /2}$.

\end{lemma}
\begin{proof}
Consider all of the $2^{|S|}$ possible assignments of sign to the entries of $\eps$ on $S$.  For each, choose an $\Splus$ and $\Sminus$ that are disjoint, have cardinality $\lceil 0.001 m \rceil$, and are supersets of the indices assigned a positive or negative sign, respectively.  Let $\Etilde_{S, x_0}$ be the event on which all sign assignments yield a $Y$ satisfying \eqref{lambda-condition},  $\|\YT\|_\text{F} \leq 1/4$, and $\| \YTp - \frac{17}{10} \ITp \| \leq 3/10$.  By Lemma \ref{lemma-dual-certificate-construction}, this event has probability at least $$1 - 2^{|S|} e^{-\gamma^*} \geq 1 - e^{-\gamma^* m /2}.$$

\end{proof}

We now show that for a fixed support set of gross errors, there is a high probability that a dual certificate exists simultaneously for all signals and for all gross errors.

\begin{lemma} \label{dual-certificate-universality-over-x}
Fix a support set $S$.  If $m \geq max(c, 4 \log(201)/\gamma^*) n$ and $|S|/m \leq \min(0.001, \gamma^* / 2 \log 2)$, then on an event of probability at least $1 - e^{-\gamma^* m/4}$, for all $x_0$ and for all $\eps$ with $\supp(\eps) \subseteq S$, there exists $Y = \A^* \lambda$ satisfying \eqref{lambda-condition} with $\alpha = 7/m$ and $\|\YT\|_\text{F} \leq 0.44$ and $\| \YTp - \frac{17}{10} \ITp \| \leq 4/10$.

\end{lemma}

\begin{proof}
By Lemma \ref{dual-certificate-over-all-errors}, for any fixed all $x_0$ such that $\|x_0\|=1$, for all $\eps$ with $\supp(\eps) \subseteq S$, there exists a $Y = \A^* \lambda$ such that
\begin{align}
\| \lambda \|_\infty &\leq \frac{7}{m} \\
\lambda_S &= \frac{7}{m} \sgn(\eps_S)\\
\| \YTp + 1.7 \ITp\| &\leq 0.3 \\
\| \YT\|_\text{F} &\leq  0.25
\end{align}
on the event $\Etilde_{S, x_0}$, which has probability at least $1 - e^{-\gamma^* m/2}$.  By Lemma 5.2 in \cite{V2012}, there exists a net $\Neps$ such that $|\Neps| \leq (1 + 2/\eps)^n$.  Hence, such a $Y$ exists simultaneously for all $x_0 \in \Neps$ on an event of probability at least $1 - (1+2/\eps)^n e^{-\gamma^* m/2}$.  If $m > 4  n\log ( 1 + 2/\eps) / \gamma^*$, then such a $Y$ exists simultaneously for all $x_0 \in \Neps$ with probability at least $1 - e^{-\gamma^* c n / 4} \geq 1 - e^{-\gamma^* m/4}$.

We now appeal to a continuity argument to show that a dual certificate exists for points not on the net $\Neps$. 
For an arbitrary $x$ such that $\|x\|_2 = 1$, we consider the $Y$ corresponding to the nearest $x_0 \in \Neps$.  Note that $\|x-x_0\| \leq \eps$ by definition of the net $\Neps$.  We now closely follow the proof  and notation of Corollary 2.4 in \cite{CL2012} to show that $Y$ is a satisfactory approximate dual certificate for $x$.

Note that $\|Y\| \leq 2.5$.  Let $\Delta  = \xxt - \xoxot$ and note that $\|\Delta\|_\text{F} \leq 2 \eps.$  Let $T = T_x$.  
Now we have 
\begin{align}
 \YTp + 1.7 \ITp &= \YTop + 1.7 \ITop - R_1 \\
 \YT &= \YTo + R_2
\end{align}
where
\begin{align}
R_1&=  \Delta Y(I - \xoxot) + (I - \xoxot)Y\Delta - \Delta Y \Delta + 1.7 \Delta\\
R_2&= \Delta Y (I - \xoxot) + (I - \xoxot) Y \Delta - \Delta Y \Delta
\end{align}

We observe that \begin{align}
\|R_1\| &\leq 2 \|Y\| \|\Delta\| \| I - \xoxot \| + \|Y\| \|\Delta\|^2 + 1.7 \|\Delta\| \leq 13.4 \eps + 10 \eps^2\\
\|R_2\| &\leq \sqrt{2} \|R_2\| \leq \sqrt{2} \left( 2 \|Y\| \|\Delta \| \|I - \xoxot\| + \|Y\| \|\Delta \|^2 \right) \leq 10 \sqrt{2} ( \eps + \eps^2)
\end{align}
If we choose $\eps = 0.01$, $\|R_1\| \leq 0.135$ and $\|R_2\|_\text{F} \leq 0.143$. and
\begin{align}
\|\YTp + 1.7 \ITp\| &\leq 0.3 + 0.135 = 0.435 \\
\| \YT\|_\text{F} &\leq 0.25 + 0.143 = 0.393
\end{align}

\end{proof}

We can now prove Theorem \ref{theorem-exact-recovery}.
\begin{proof}[Proof of Theorem \ref{theorem-exact-recovery}]
Assume that $|S|/m \leq \min(0.001, \gamma^*/2 \log 2)$ and $m \geq \max(c, 4 \log(201)/\gamma^*)$. By Lemma \ref{dual-certificate-universality-over-x}, there is an event of probability at least $1 - e^{-\gamma^* m/4}$ such that for all $x_0$ and for all $\eps$ with $\supp(\eps) \subseteq S$, there exists a $Y = \A^* \lambda$ satisfying \eqref{ytp-condition}--\eqref{lambda-condition}.  Choose a superset $\overline{S}$ such that $|\overline{S}| = 0.01m$.  On the intersection of this event with $E_{\overline{S}}$, Lemma \ref{lemma-recovery-with-dual-certificate} guarantees that $\| \Xhat - X_0 \|_\text{F} \leq C \|\eta\|_1 / m$.  The intersection of these events has probability at least $ 1 - e^{-\gamma m}$ for some $\gamma$. 

The proof of the error estimate $\|\hat{x} - x_0 \| \leq C' \min \left( \|x_0\|, \frac{\|\eta\|_1}{m \|x_0\|} \right)$ can be found in \cite{CSV2013}.
\end{proof}

\bibliographystyle{plain}
\bibliography{arberror-refs}

\end{document}